\documentclass{amsart}
\usepackage{amsfonts}

\usepackage{graphicx}
\usepackage{amscd}

\newtheorem{theorem}{Theorem}
\theoremstyle{plain}

\newtheorem{conclusion}{Conclusion}

\numberwithin{equation}{section}

\begin{document}
\title[Second order condition and quantile representations]{Another look at Second order condition in Extreme Value Theory}
\author[G. S. LO]{Gane Samb LO$^{*}$}
\address{$^{*}$ LSTA, UPMC, France and LERSTAD, Universit\'e Gaston Berger de
Saint-Louis, SENEGAL}
\email{gane-samb.lo@ugb.edu.sn, ganesamblo@ufrsat.org}
\urladdr{www.lsta.upmc.fr}
\author{Adja Mbarka FALL$^{**}$}
\address{$^{**}$ MAPMO, Universit\'e d'Orl\'eans, France and LERSTAD, Universit\'e
Gaston Berger de Saint-Louis, SENEGAL}
\email{amfall@ufrsat.org}

\begin{abstract}
This note compares two approaches both alternatively used when establishing
normality theorems in univariate Extreme Value Theory. When the underlying
distribution function ($df$) is the extremal domain of attraction, it is
possible to use representations for the quantile function and regularity
conditions (RC), based on these representations, under which strong and weak
convergence are valid. It is also possible to use the now fashion second
order condition (SOC), whenever it holds, to do the same. Some authors
usually favor the first approach (the SOC one) while others are fond of the
second approach that we denote as the representational one. This note aims
at comparing the two approaches and show how to get from one to the other.
The auxiliary functions used in each approach are computed and compared.
Statistical applications using simultaneously both approaches are provided.
A final comparison is provided.
\end{abstract}

\maketitle

\section{Introduction}

\noindent Statistical modelling based on the univariate Extreme Values
Theory usually requires regularity conditions of the underlying
distributions. Since the work of de Haan and Stadtm\"{u}ller (\cite
{haanstadm} ) on the the so-called Second Order Condition (SOC), using this
SOC has became fashion in research papers so that this way of doing is the
mainstream one, led by de Haan (see for instance \ \cite{deh}, \cite{segers}%
, \cite{hf}). However the second order condition does not always hold as we
will show it (see (\ref{exp1})), although a large class of distribution
functions fulfills it. Yet, there exists an other approach, that is the
representational one, based on the Karamata representation for a slowly
varying function. In this view, each distribution function $F$ in the
extremal domain may be represented a couple of functions $p(s)$ and $b(s)$, $%
s\in (0,1)$, to be precised in Theorem \ref{A}. This approach is the one
preferred by many other authors, for instance Cs\"{o}rg\H{o}, Deheuvels, and
Mason, (\cite{cdm}), Lo (\cite{gslo}), Hall (see \cite{hall1},\cite{hall2})
etc. This latter in particularly adapted for the use on the Gaussian
approximations like that of Cs\"{o}rg\H{o}-Cs\"{o}rg\H{o}-Horv\`{a}th-Mason (
\cite{cchm}).\medskip

\noindent This motivates us to undertake here a comparative study of the
second order condition in the two approaches and provide relations and
methods for moving from one to the other. We give specific statistical
applications \ using simultanuously the two ways. The paper is to serve as a
tool for comparitive reading of papers based on the two approaches.\medskip

\noindent The paper is organized as follows. In Section 1, we introduce the
second order condition in the frame of de Haan and Stadm\"{u}ller (\cite
{haanstadm} ) using quantile functions. In Section 3, we recall the
representational scheme and link them to the second order condition.
Precisely, we express the second order condition, when it holds, through the
couple of functions $(p,b)$ associated with a $df$ attracted to the extremal
domain. The results are then given through the $df$ $G(x)=F(e^{x})$, $x\in R$%
, that is the most used in statistical context. In Section 4, we settle a
new writing the SOC for the quantiles while the auxiliary functions of that
condition, denoted as $s$ and $S$, are computed for a large number of $%
df^{\prime }s$. In Section 6, we deal with applications in statistical
contexts. The first concerns the asymptotic normality of the large quantiles
process and the second treats the functional Hill process. In both cases, we
use the two approaches. We finish by comparing the two methods at the light
of these applications.

\section{\textbf{The second order condition}}

\subsection{Definition and expressions.}

\noindent Consider a $df$ $F$ lying in the extremal domain of attraction of
the Generalized Extreme Value (GEV) distribution, that is 
\begin{equation*}
G_{\gamma }(x)=exp(-(1+\gamma x)^{-1/\gamma }),\newline
for\newline
1+\gamma x>0,
\end{equation*}

\noindent denoted denoted $F\in D(G_{\gamma })$, and let $%
U(x)=(1/(1-F))^{-1} $, where for any nondecreasing and right-continuous
function $L:R\mapsto \lbrack a,b]$, with $a<b$, 
\begin{equation*}
L^{-1}(t)=\inf \{x\geq t,L(x)\geq t\},a\leq t\leq b,
\end{equation*}
is the generalized inverse of $L$. One proves (see \cite{hf}, p. 43) that
there exists a positive function $a(t)$ of $t\in R$, such that 
\begin{equation*}
\forall (x>0),\lim_{t\rightarrow \infty }\frac{U(xt)-U(t)}{a(t)}=\frac{%
x^{\gamma }-1}{\gamma }=D_{\gamma }(x),
\end{equation*}

\noindent where $(x^{\gamma }-1)/\gamma $ is interpreted as $\log (x)$ for $%
\gamma =0.$ Now, by definition, $F$ is said to satisfy a second order
condition (\cite{haanstadm}) if and only if there exists a function $A(t)$
of $t\in R$ with a constant sign such that 
\begin{equation}
\lim_{t\rightarrow \infty }\frac{\frac{U(xt)-U(t)}{a(t)}-\frac{x^{\gamma }-1%
}{\gamma }}{A(t)}=H(x),  \tag{SOCU}
\end{equation}
holds. According to Theorem 2.3.3 in (de haan and Ferreira), the function $H$%
, when it is not a multiple of $D_{\gamma }(x),$ can be written as 
\begin{equation*}
H_{\gamma ,\rho }(x)=c_{1}\int_{1}^{x}s^{\gamma -1}\int_{1}^{s}u^{\gamma -1}%
\text{ }du\text{ }ds+c_{2}\int_{1}^{x}s^{\gamma +\rho -1}ds,
\end{equation*}
where $\rho $ is a negative number and the functions $a(t)$ and $A(t)$
satisfies for any $x>0$, 
\begin{equation*}
\lim_{t\rightarrow \infty }\frac{a(tx)/a(t)-x^{\gamma }}{A(t)}%
=c_{1}x^{\gamma }\frac{x^{\rho }-1}{\rho }
\end{equation*}
and 
\begin{equation*}
\lim_{t\rightarrow \infty }A(tx)/A(t)=x^{\rho }.
\end{equation*}
According to Corollary 2.3.5 in de Haan and Ferreira, one can choose a
positive function $a^{\ast }(t)$ and a function $A^{\ast }(t)$ of constant
sign such that 
\begin{equation*}
\frac{\frac{U(xt)-U(t)}{a^{\ast }(t)}-\frac{x^{\gamma }-1}{\gamma }}{A^{\ast
}(t)}\rightarrow H_{\gamma ,\rho }^{\ast }(x),
\end{equation*}
with 
\begin{equation*}
H_{\gamma ,\rho }^{\ast }(x)=\left\{ 
\begin{array}{ll}
\frac{x^{\gamma +\rho }-1}{\gamma +\rho },\text{ } & \gamma +\rho \neq
0,\rho <0, \\ 
\log x,\text{ \ } & \gamma +\rho =0,\rho <0, \\ 
\frac{1}{\gamma }x^{\gamma }\log x,\text{ \ } & \gamma \neq 0,\rho =0, \\ 
(\log (x))^{2}/2,\text{ \ \ \ \ \ } & \gamma =\rho =0,
\end{array}
\right.
\end{equation*}
\begin{equation*}
a^{\ast }(t)=\left\{ 
\begin{array}{lc}
a(t)(1-A(t)/\rho ), & \rho >0, \\ 
a(t)(1-A(t)/\gamma ),\text{ \ } & \rho =0,\gamma \neq 0, \\ 
a(t),\text{ \ \ \ \ } & \gamma =\rho
\end{array}
\right.
\end{equation*}
and 
\begin{equation*}
A^{\ast }(t)=A(t)/\rho I(\rho >0)+A(t)I(\rho =0).
\end{equation*}

\noindent To see that the $SOC$ does not necessary hold, consider the
standard exponential distribution function. We have $U(t)=log(t)$ and 
\begin{equation}
U(tx)-U(t)=\log x=[(x^{\gamma }-1)/\gamma ]_{\gamma =0}.  \label{exp1}
\end{equation}
It is clear that the function $a$ is necessarily constant and equal to the
unity and the second order condition is here meaningless. As a consequence,
the results obtained under a second order condition are partial.

\subsection{Expression in terms of generalized inverse functions.}

\noindent We are going to express the SOC through the generalized function $%
F^{-1}(1-u).$ Let 
\begin{equation*}
F\in D(G_{\gamma }).
\end{equation*}
With this parameterization, the case $\gamma =0$ corresponds to $D(\Lambda )$%
, the case $-\infty <\gamma <0$ to $D(\psi _{1/\gamma })$ and finally, the
case $0<\gamma <+\infty $ to $D(\phi _{1/\gamma }).$\ The second order
condition will become : there exist a positive function $s(u)$ and a
function $S(u)$ with constant sign such that for any $x>0,$%
\begin{equation}
\frac{\dfrac{F^{-1}(1-ux)-F^{-1}(1-u)}{s(u)}-\gamma ^{-1}(x^{\gamma }-1)}{%
S(u)}=H_{\gamma ,\rho }(1/x)=h_{\gamma ,\rho }(x).  \tag{SOCF}
\end{equation}

\section{Representation for F$\in D(G_{\protect\gamma })$}

\subsection{Representations}

\noindent Now we recall the classical representations of $df$ attracted to
some nondegenerated extremal $df.$

\begin{theorem}
\label{A} We have :

\begin{enumerate}
\item  Karamata's representation (KARARE)

\begin{enumerate}
\item[(a)]  If $F\in D(G_{\gamma }),$ $\gamma >0$, then there exist two
measurable functions $p(u)$ and $b(u)$ of $u\in (0,1)$ such that $\sup
(\left\vert p(u)\right\vert ,\left\vert b(u)\right\vert )\rightarrow 0$ as $%
u\rightarrow 0$ and a positive constant c so that 
\begin{equation}
G^{-1}(1-u)=\log c+\log (1+p(u))-\gamma \log u+(\int_{u}^{1}b(t)t^{-1}dt),%
\text{ }0<u<1,  \label{rep1}
\end{equation}
where $G^{-1}(u)=\inf \{x,G(x)\geq u\},$ $0<u\leq 1$ is the generalized
inverse of $G$ with $G^{-1}(0)=G^{-1}(0+)$.\medskip

\item[(b)]  If $F\in D(G_{\gamma }),$ $\gamma <0$, then $y_{0}(G)=\sup \{x,$ 
$G(x)<1\}<+\infty $ and there exist two measurable functions $p(u)$ and $b(u)
$ for $u\in (0,1)$ and a positive constant c as defined in (\ref{rep1}) such
that 
\begin{equation}
y_{0}-G^{-1}(1-u)=c(1+p(u))u^{-\gamma }\exp (\int_{u}^{1}b(t)t^{-1}dt),0<u<1.
\label{rep2}
\end{equation}
\end{enumerate}

\item  Representation of de Haan (Theorem 2.4.1 in \cite{dehaan}),

\noindent If $G\in D(G_{0})$, then there exist two measurable functions $%
p(u) $ and $b(u)$ of $u\in (0,1)$ and a positive constant c as defined in (%
\ref{rep1}) such that for 
\begin{equation}
s(u)=c(1+p(u))\exp (\int_{u}^{1}b(t)t^{-1}dt),\text{ }0<u<1,  \label{rep3b}
\end{equation}
we have for some constant $d\in R$,

\begin{equation}
G^{-1}(1-u)=d-s(u)+\int_{u}^{1}s(t)t^{-1}dt,0<u<1.  \label{rep3a}
\end{equation}
\end{enumerate}
\end{theorem}

\noindent It is important to remark at once that any $df$ in the extremal
domain of attraction is associated with a couple of functions $(p,b)$ used
in each appropriate representation.

\subsection{Preparation of second order condition.}

\noindent We are now proving, under $F\in D(G_{\gamma })$ that for $x>0$, 
\begin{equation*}
\lim_{u\rightarrow \infty }\frac{F^{-1}(1-ux)-F^{-1}(1-u)}{s(u)}=d_{\gamma
}(x)=\gamma ^{-1}(x^{-\gamma }-1).
\end{equation*}

\subsubsection{F $\in D(G_{0})$}

\noindent The representation (\ref{rep3a}) is valid for $F^{-1}(1-u)$, $u\in
(0,1)$. We get for $u\in (0,1)$ and $x>0$ : 
\begin{equation*}
F^{-1}(1-xu)-F^{-1}(1-u)=s(u)-s(ux)+\int_{ux}^{u}s(t)/t\text{ }dt.
\end{equation*}
For $v\in \lbrack \min (ux,u),\max (ux,u)]=A(u,x),$%
\begin{equation*}
s(v)/s(u)=\frac{1+p(t)}{1+p(u)}\exp (\int_{v}^{u}b(t)/t\text{ }dt).
\end{equation*}
By letting $pr(u,x)=\sup \{\left\vert p(t)\right\vert ,0\leq t\leq \max
(ux,u)\}$ and $br(u,x)=\sup \{\left\vert b(t)\right\vert ,0\leq t\leq \max
(ux,u)\},$ one quickly shows that, for u sufficiently small, 
\begin{equation*}
\sup_{t\in A(u,x)}\left\vert 1-(1+p(v))/(1+p(u)\right\vert \leq 2pr(u,x).
\end{equation*}
and 
\begin{equation*}
x^{-br(u,x)}\leq \exp (\int_{ux}^{u}b(t)/t\text{ }dt)\leq x^{br(u,x)}
\end{equation*}
and then 
\begin{equation*}
\sup_{v\in A(u,x)}\left\vert 1-\exp (\int_{v}^{u}b(t)/t\text{ }%
dt)\right\vert =O(-br(u,x)).
\end{equation*}
\begin{equation*}
\leq \left\vert 1+x\right\vert (1\vee \left\vert x\right\vert )^{1+br(u,x)}
\end{equation*}
It follows that 
\begin{equation*}
\sup_{v\in A(u,x)}\left\vert 1-s(v)/s(u)\right\vert =O(\max
(pr(u,x),br(u,x))\rightarrow 0,
\end{equation*}
as $u\rightarrow 0$. We get, for $pbr(u,x)=pr(u,x)\times br(u,x),$%
\begin{equation*}
\frac{F^{-1}(1-xu)-F^{-1}(1-u)}{s(u)}+\log x
\end{equation*}
\begin{equation*}
=O(pbr(u,x))+O(pbr(u,x)\log x)\rightarrow 0.
\end{equation*}
Then 
\begin{equation*}
\lim_{u\rightarrow 0}\frac{F^{-1}(1-xu)-F^{-1}(1-u)}{s(u)}=-\log x=\left[
\gamma (x^{-1/\gamma }-1)\right] _{\gamma =\infty }.
\end{equation*}
We notice that 
\begin{equation}
\frac{F^{-1}(1-xu)-F^{-1}(1-u)}{s(u)}+\log x  \label{gum01}
\end{equation}
\begin{equation*}
=O(p(u,x)+b(u,x)),
\end{equation*}
where 
\begin{equation}
p(u,x)=1-(1-p(ux))/(1-p(u))\rightarrow 0\text{ }as\text{ }u\rightarrow 0
\label{gum02}
\end{equation}
and 
\begin{equation}
b(u,x)=\int_{ux}^{u}\frac{1}{t}\left[ \frac{1-p(t)}{(1+p(u)}\exp
(\int_{u}^{t}v^{-1}b(v)\text{ }dv)-1\right] dt\rightarrow 0\text{ }as\text{ }%
u\rightarrow 0.  \label{gum03}
\end{equation}

\subsubsection{F$\in D(G_{1/\protect\gamma })$ $(\protect\gamma >0)$}

\noindent We have the KARARE representation 
\begin{equation*}
F^{-1}(1-u)=c(1+p(u))\text{ }u^{-\gamma }\text{ }exp(\int_{u}^{1}t^{-1}b(t)%
\text{ }dt),\text{{}}u\in (0,1).
\end{equation*}
Then, for $s(u)=\gamma cu^{-\gamma }$ $exp(\int_{u}^{1}t^{-1}b(t)$ $dt)$,
for $x>0$ 
\begin{equation*}
\frac{F^{-1}(1-ux)-F^{-1}(1-u)}{s(u)}=\gamma ^{-1}\left\{
(1+p(ux))x^{-1/\gamma }\text{ }exp(\int_{ux}^{x}t^{-1}b(t)\text{ }%
dt)-1-p(u)\right\} .
\end{equation*}
As previously, we readily see that 
\begin{equation*}
exp(\int_{ux}^{x}t^{-1}b(t)\text{ }dt)=\exp (O(br(u,x)\log
x))=1+O(br(u,x)\log x)\rightarrow 1,
\end{equation*}
as $u\rightarrow 0.$ It follows that 
\begin{equation}
\frac{F^{-1}(1-ux)-F^{-1}(1-u)}{s(u)}\rightarrow \gamma ^{-1}(x^{-\gamma
}-1).  \label{par01}
\end{equation}
Moreover we have 
\begin{equation*}
\frac{F^{-1}(1-ux)-F^{-1}(1-u)}{s(u)}-\gamma ^{-1}(x^{-\gamma }-1)
\end{equation*}
\begin{equation}
=\gamma ^{-1}\left\{ x^{\gamma }\left( (1+p(ux))\text{ }exp(%
\int_{ux}^{x}t^{-1}b(t)\text{ }dt)-1\right) -p(u)\right\} =pb(u,x).
\label{par02}
\end{equation}
Notice that we may also take $s(u)=F^{-1}(1-u).$

\subsubsection{F$\in D(G_{\protect\gamma }),\protect\gamma <0$}

\noindent We have $x_{0}=\sup \{x,F(x)<1\}<+\infty $ and the following
representation holds : 
\begin{equation*}
x_{0}-F^{-1}(1-u)=c(1+p(u))u^{-\gamma }\exp (\int_{u}^{1}t^{-1}b(t)dt),\text{%
{}}u\in (0,1).
\end{equation*}
Then, for $s(u)=-\gamma cu^{-\gamma }\exp (\int_{u}^{1}t^{-1}b(t)dt)$, $u\in
(0,1)$ and $x>0$, 
\begin{equation*}
\frac{F^{-1}(1-ux)-F^{-1}(1-u)}{s(u)}=\frac{%
(x_{0}-F^{-1}(1-u))-(x_{0}-F^{-1}(1-ux))}{s(u)}
\end{equation*}
\begin{equation*}
=-\gamma ^{-1}\left\{ 1+p(u)-x^{-\gamma }(1+p(ux))\exp
(\int_{ux}^{x}t^{-1}b(t)dt)\right\}
\end{equation*}
\begin{equation*}
=\gamma ^{-1}\left\{ x^{-\gamma }(1+p(ux))\exp
(\int_{ux}^{x}t^{-1}b(t)dt)\right\} -1-p(u)
\end{equation*}
\begin{equation*}
\rightarrow \gamma ^{-1}(x^{-\gamma }-1),
\end{equation*}
as $u\rightarrow 0$. Likely to the case $0<\gamma <\infty ,$%
\begin{equation}
\frac{F^{-1}(1-ux)-F^{-1}(1-u)}{s(u)}-\gamma (x^{-1/\gamma }-1)
\label{wei01}
\end{equation}
\begin{equation}
=\gamma ^{-1}\left\{ x^{-\gamma }\left( (1+p(ux))\text{ }exp(%
\int_{ux}^{x}t^{-1}b(t)\text{ }dt)-1\right) -p(u)\right\} =pb(u,x).
\label{wei02}
\end{equation}
Notice that we may also take $s(u)=x_{0}-F^{-1}(1-u).$

\section{Second order condition via representations}

\subsection{Case by case}

\subsubsection{$F\in D(G_{\protect\gamma }),0<\protect\gamma <\infty $}

\noindent The second order condition is equivalent to finding a fonction $%
S(u)$ of constant sign such that 
\begin{equation*}
S(u)^{-1}\left\{ \frac{F^{-1}(1-ux)-F^{-1}(1-u)}{s(u)}-\gamma
^{-1}(x^{-\gamma }-1)\right\} =S(u)^{-1}pb(u,x)
\end{equation*}
converges to a function $h_{\gamma ,\rho }$, where 
\begin{equation*}
pb(u,x)=\gamma \left\{ x^{-\gamma }\left( (1+p(ux))\text{ }%
exp(\int_{ux}^{x}t^{-1}b(t)\text{ }dt)-1\right) -p(u)\right\}
\end{equation*}

\noindent and $s(u)$ may be taken $as$ $F^{-1}(1-u).$

\subsubsection{$F\in D(G_{\protect\gamma }),-\infty <\protect\gamma <0$}

\noindent The second order condition is equivalent to finding a fonction $%
S(u)$ of constant sign such that 
\begin{equation*}
S(u)^{-1}\left\{ \frac{F^{-1}(1-ux)-F^{-1}(1-u)}{s(u)}-\gamma
^{-1}(x^{-\gamma }-1)\right\} =S(u)^{-1}pb(u,x)
\end{equation*}
converges to a function $h_{\gamma ,\rho }$ where 
\begin{equation*}
pb(u,x)=\gamma ^{-1}\left\{ x^{-\gamma }\left( (1+p(ux))\text{ }%
exp(\int_{ux}^{x}t^{-1}b(t)\text{ }dt)-1\right) -p(u)\right\}
\end{equation*}

\noindent and $s(u)$ may be taken $as$ $x_{0}-F^{-1}(1-u)$ and $x_{0}$ is
the upper endpoint of $F$.

\subsubsection{$F\in D(G_{0})$}

\noindent The second order condition is equivalent to finding a fonction $%
S(u)$ of constant sign such that 
\begin{equation*}
S(u)^{-1}\left\{ \frac{F^{-1}(1-ux)-F^{-1}(1-u)}{s(u)}-\log x)\right\}
=S(u)^{-1}(p(u,x)+b(u,x))
\end{equation*}
converges to a function $h_{1/\gamma ,\rho },$ where

\begin{equation*}
p(u,x)=1-(1-p(ux))/(1-p(u))\rightarrow 0\text{ }as\text{ }u\rightarrow 0
\end{equation*}
and 
\begin{equation*}
b(u,x)=\int_{ux}^{u}\frac{1}{t}(\frac{1-p(t)}{(1+p(u)}\exp
(\int_{u}^{t}v^{-1}b(v)\text{ }dv)-1)dt\rightarrow 0\text{ }as\text{ }%
u\rightarrow 0.
\end{equation*}
and $s(u)$ may be taken $u^{-1}\int_{u}^{1}(1-s)dF^{-1}(s).$

\section{Special cases}

\subsection{Statistical context.}

\noindent In the statistical context, especially in the exteme value index
estimation, the bulk of the work is done with 
\begin{equation*}
G^{-1}(1-u)=\log F^{-1}(1-u).
\end{equation*}
Let $F\in D(G_{\gamma })$. The three cases $-\infty <\gamma <0,$ $0<\gamma
<+\infty ,$ $\gamma =0$ respectively imply 
\begin{equation*}
G\in D(G_{1/\gamma })\text{,}
\end{equation*}
\begin{equation*}
G\in D(G_{0})\text{ and }s(u,G)\rightarrow \gamma
\end{equation*}
and 
\begin{equation*}
G\in D(G_{0})\text{ and }s(u,G)\rightarrow 0.
\end{equation*}
For $F\in D(G_{\gamma }),$ $0<\gamma <\infty ,$ we have a representaion like 
\begin{equation*}
G^{-1}(1-u)=c+\log (1+p(u))-\gamma \log u+\int_{u}^{1}t^{-1}b(t)dt.
\end{equation*}
We take here $s(u)=\gamma $. The second order conditions becomes 
\begin{eqnarray*}
S(u)^{-1}\left\{ \frac{G^{-1}(1-ux)-G^{-1}(1-u)}{\gamma }-\log x\right\}
&=&S(u)^{-1}\left\{ \gamma ^{-1}\log \frac{1+p(ux)}{1+p(u)}+\gamma
^{-1}\int_{ux}^{u}t^{-1}b(t)dt\right\} \\
&=&A(u)^{-1}\text{ }pb(u,x)\rightarrow h_{0,\rho }(x).
\end{eqnarray*}
Denote $dG^{-1}(1-u)/du=G^{-1}(1-u)^{\prime }$ whenever if exists. Now if $%
G^{-1}(1-u)^{\prime }$ exists for $u$ near zero, we may take 
\begin{equation*}
b(u)=G^{-1}(1-u)^{\prime }+\gamma \rightarrow 0,
\end{equation*}
For $F\in D(G_{\gamma }),$ $-\infty <\gamma <0,$ we may transfer the SOC to $%
G$ in a way similar as to $F$, with 
\begin{equation*}
\log x_{0}-G^{-1}(1-u)=c(1+p(u))u^{\gamma }\exp (\int_{u}^{1}t^{-1}b(t)dt).
\end{equation*}

\noindent For $F\in D(G_{0}).$ If $s(u)=u(G^{-1}(1-u)^{\prime })\rightarrow
0,$ we will have 
\begin{equation*}
G^{-1}(1-s)=d-\int_{u}^{1}t^{-1}s(t)dt.
\end{equation*}
We may take
\begin{equation*}
b(u)=us^{\prime }(u).
\end{equation*}
The second order condition becomes simpler as 
\begin{eqnarray*}
S(u)^{-1}b(u,x) &=&S(u)^{-1}\int_{ux}^{u}\frac{1}{t}(\frac{1-p(t)}{(1+p(u)}%
\exp (\int_{u}^{t}v^{-1}b(v)\text{ }dv)-1)dt \\
&\rightarrow &h_{0,\rho }(x).
\end{eqnarray*}
Moreover, for $g(x)=dG(x)/dx,$ if 
\begin{equation*}
b(u)=us^{\prime }(u)/s(u)=1-u(G^{-1}(1-u)g(G^{-1}(1-u))^{-1}\rightarrow 0,
\end{equation*}
we have 
\begin{equation*}
s(u)=c\exp (\int_{u}^{1}t^{-1}b(t)dt),
\end{equation*}
and the SOC becomes 
\begin{eqnarray*}
S(u)^{-1}b(u,x) &=&S(u)^{-1}\int_{ux}^{u}\frac{1}{t}\left[ \exp
(\int_{u}^{t}(\frac{1}{v}-\frac{1}{G^{-1}(1-\nu )g(G^{-1}(1-\nu )})-1)d\nu %
\right] dt \\
&\rightarrow &h_{0,\rho }(x).
\end{eqnarray*}

\section{Finding the functions $b$ and $S.$}

\subsection{Determination of the function $b$}

\noindent In the usual cases, the function is ultimately differentiable,
that is in a right neigbourhood of $x_{0}(F)$. It is then easy to find the
function \ $b$ by derivating $G^{-1}(1-u).$ In summary, for $D(G_{\gamma }),$
$\gamma >0,$ the function $b$ in the representation of $G^{-1}(1-u)$ is 
\begin{equation*}
b(u)=-\left\{ \gamma +u(G^{-1}(1-u))^{\prime }\right\} .
\end{equation*}
The function $b$ in the representation of $F^{-1}(1-u)$ is defined by
\begin{equation*}
-b(u)=\gamma +u(G^{-1}(1-u))^{\prime }/G^{-1}(1-u))
\end{equation*}
For $\gamma =+\infty ,$ the function $b$ is the representation of $%
G^{-1}(1-u)$, that is, 
\begin{equation*}
b(u)=-us^{\prime }(u)/s(u),
\end{equation*}
where
\begin{equation*}
s(u)=u(G^{-1}(1-u))^{\prime }.
\end{equation*}
For $\gamma <0$ and $y_{0}(G)=y_{0},$%
\begin{equation*}
b(u)=-\gamma +u(G^{-1}(1-u))^{\prime }/(x_{0}-G^{-1}(1-u))
\end{equation*}

\noindent Then we apply these formulas and determine the function $b$ for
usual $df$'s. Regularity conditions in the representational approach mainly
rely on the function $b$, while they rely of the function $S$ for the SOC
approach. It is then interesting to have both functions for usual $df$'s in
tables in Subsection \ref{sstable}, following \cite{segers}.

\bigskip

\subsection{The function $S$ for the second order condition}

\noindent Functions $a$ and $A$ in the $(SOCU)$, as well as the functions $%
H_{\gamma ,\rho }$ are available in the usual cases (see \cite{segers} for
example). It is not the case for the $(SOCF)$ expressed in terms of the
quantile functions. We then seize this opportunity to compute their analogs $%
s$ and $S$ in the this case for the usual $df$'s. The results are summarized
in our tables in Subsection \ref{sstable}.

\subsubsection{The Singh-Maddala Law}

\noindent Let for constants $a,$ $b$ and $c$, for $x\geq 0,$%
\begin{equation*}
1-F(x)=(1+ax^{b})^{-c},
\end{equation*}
the so-called Singh-Madalla $df$. This function plays a special role in
income fitting distribution. It is clear that 
\begin{equation*}
F\in D(G_{bc}).
\end{equation*}
Put 
\begin{equation*}
\lambda =1/bc.
\end{equation*}
Straightforward calculations give 
\begin{eqnarray*}
G^{-1}(1-u) &=&-b^{-1}\log a-\gamma \log u+b^{-1}\log (1-u^{1/c}) \\
&=&d-\gamma \log u-\frac{1}{b}\log (1-u_{0})+\int_{u}^{u_{0}}t^{-1}B(t)dt,
\end{eqnarray*}
where 
\begin{equation*}
a=-b^{-1}\log a,
\end{equation*}

\begin{equation*}
B(u)=\gamma u^{1/c}(1-u^{1/c})^{-1}
\end{equation*}
and $u_{0}\in ]0,1[.$ Put $K_{0}=-\frac{1}{b}\log (1-u_{0})/u_{0}$ and 
\begin{equation*}
b(u)=B(u)\mathbb{I}_{(0\leq u\leq u_{0})}+K_{0}\mathbb{I}_{(u_{0}\leq u\leq
1)},
\end{equation*}
we get 
\begin{equation*}
G^{-1}(1-u)=d-\gamma \log u+\int_{u}^{1}t^{-1}b(t)dt,
\end{equation*}
with 
\begin{equation*}
b(u)\rightarrow 0.
\end{equation*}
We have 
\begin{eqnarray*}
\frac{G^{-1}(1-ux)-G^{-1}(1-u)}{\gamma }+\log x &=&\frac{1}{\gamma b}(\log
(1-(ux)^{1/c})-\log (1-u^{1/c})) \\
&=&\frac{1}{\gamma b}(-x^{1/c}u^{1/c}+u^{1/c}+O(u^{2/c})).
\end{eqnarray*}
Thus, for $S(u)=cu^{1/c}/(\gamma b),$ we get 
\begin{equation*}
\frac{\dfrac{G^{-1}(1-ux)-G^{-1}(1-u)}{\gamma }-\log x}{S(u)}=\frac{x^{1/c}-1%
}{-1/c}=h_{0,\rho }(x)=H_{0,\rho }(1/x).
\end{equation*}
This corresponds to a second order condition. As for $F^{-1}(1-u)$ itself,
we have
\begin{equation*}
F^{-1}(1-u)=(1/a)^{1/b}(1-u^{-1/c})^{1/b}
\end{equation*}
and 
\begin{equation*}
b(u)=F^{-1}(1-u)^{\prime }/F^{-1}(1-u)+\frac{1}{bc}=(1-\frac{u^{-1/c}}{%
1+u^{-1/c}})/bc\rightarrow 0.
\end{equation*}

\noindent We have for $s(u)=F^{-1}(1-u)/(bc),$%
\begin{equation*}
\frac{F^{-1}(1-ux)-F^{-1}(1-u)}{s(u)}=bc(\left( \frac{1-x^{-1/c}u^{-1/c}}{%
1-u^{-1/c}}\right) ^{1/b}-1)\rightarrow \frac{x^{-1/bc}-1}{1/bc}.
\end{equation*}
Next
\begin{equation*}
\frac{F^{-1}(1-ux)-F^{-1}(1-u)}{s(u)}-\frac{x^{-1/bc}-1}{1/bc}=(bc)\left\{
\left( \frac{1-x^{-1/c}u^{-1/c}}{1-u^{-1/c}}\right) ^{1/b}-\left(
x^{-1/c}\right) ^{1/b}\right\}
\end{equation*}
\begin{equation*}
=(bc)b^{-1}\left\{ \frac{1-x^{-1/c}u^{-1/c}}{1-u^{-1/c}}-\frac{%
x^{-1/c}(1-u^{-1/c})}{1-u^{-1/c}}\}\right\} \zeta (u,x)^{(1/b)-1},
\end{equation*}
where $\zeta (u,x)\in I(x^{-1/c},b),$ with $%
b=(1-x^{-1/c}u^{-1/c})/(1-u^{-1/c})\sim x^{-1/c}.$ Hence
\begin{eqnarray*}
\frac{F^{-1}(1-ux)-F^{-1}(1-u)}{s(u)}-\frac{x^{-1/bc}-1}{1/bc} &=&c\left\{ 
\frac{1-x^{-1/c}u^{-1/c}}{1-u^{-1/c}}-\frac{x^{-1/c}(1-u^{-1/c})}{1-u^{-1/c}}%
\right\} \zeta (u,x)^{(1/b)-1} \\
&=&c\left\{ \frac{1-x^{-1/c}}{1-u^{-1/c}}\right\} \zeta (u,x)^{(1/b)-1}.
\end{eqnarray*}
Finally for $S(u)=-(1-u^{-1/c})^{-1},$%
\begin{equation*}
\frac{\frac{F^{-1}(1-ux)-F^{-1}(1-u)}{s(u)}-\frac{x^{-1/bc}-1}{1/bc}}{S(u)}%
\rightarrow h(x)=\frac{\left( x^{-1/c}-1\right) x^{-1/(bc)+1/c}}{1/c}
\end{equation*}

\subsubsection{Burr's df}

$1-F(x)=(x^{-\rho /\gamma }+1)^{1/\rho },\rho <0,$ and
\begin{equation*}
F^{-1}(1-u)=(u^{\rho }-1)^{-\gamma /\rho }.
\end{equation*}
For $s(u)=\gamma F^{-1}(1-u)$, 
\begin{eqnarray*}
\frac{F^{-1}(1-ux)-F^{-1}(1-u)}{s(u)} &=&\frac{1}{\gamma }\left\{ \frac{%
(x^{\rho }u^{\rho }-1)^{-\gamma /\rho }}{(u^{\rho }-1)^{-\gamma /\rho }}%
-1\right\}  \\
&=&\frac{1}{\gamma }\left\{ \left( \frac{x^{\rho }u^{\rho }-1}{u^{\rho }-1}%
\right) ^{-\gamma /\rho }-1\right\}  \\
\  &\rightarrow &\frac{x^{-\gamma }-1}{\gamma }.
\end{eqnarray*}
Next
\begin{eqnarray*}
\frac{F^{-1}(1-ux)-F^{-1}(1-u)}{s(u)}-\frac{x^{-\gamma }-1}{\gamma } &=&%
\frac{1}{\gamma }\left\{ \left( \frac{x^{\rho }u^{\rho }-1}{u^{\rho }-1}%
\right) ^{-\gamma /\rho }-x^{-\gamma }\right\}  \\
&=&\frac{1}{\gamma }\left\{ \left( \frac{x^{\rho }u^{\rho }-1}{u^{\rho }-1}%
\right) -(x^{\rho })^{-\gamma /\rho }\right\}  \\
&=&-\frac{1}{\rho }\left\{ \frac{x^{\rho }u^{\rho }-1-x^{\rho }(u^{\rho }-1)%
}{u^{\rho }-1}\right\} \zeta (u,x)^{-\gamma /\rho -1} \\
(\text{ a v\'{e}rifier})\ast \ast \ast  &=&-\frac{1}{\rho }\left\{ \frac{%
x^{\rho }u^{\rho }-1-x^{\rho }(u^{\rho }-1)}{u^{\rho }-1}\right\} \zeta
(u,x)^{-\gamma /\rho -1},
\end{eqnarray*}
where $\zeta (u,x)=I(a,b)=[a\wedge b,a\vee b],$ with $a=x^{\rho }$ and $%
b=(x^{\rho }u^{\rho }-1)/(u^{\rho }-1)\sim $ $x^{\rho }.$ Hence
\begin{equation*}
\frac{F^{-1}(1-ux)-F^{-1}(1-u)}{s(u)}-\frac{x^{-\gamma }-1}{\gamma }=-\frac{1%
}{\rho }\left\{ \frac{x^{\rho }-1}{u^{\rho }-1}\right\} \zeta (u,x)^{-\gamma
/\rho -1}.
\end{equation*}
Hence for $S(u)=(u^{\rho }-1)^{-1},$%
\begin{equation*}
\frac{\dfrac{F^{-1}(1-ux)-F^{-1}(1-u)}{s(u)}-\dfrac{x^{-\gamma }-1}{\gamma }%
}{S(u)}\rightarrow -\frac{(x^{\rho }-1)x^{-\gamma -\rho }}{\rho }.
\end{equation*}

\subsubsection{Log Exponential law}

\noindent $F(x)=1-exp(-e^{-x}),$ that is
\begin{equation*}
F^{-1}(1-u)=\log \log (1/u).
\end{equation*}
We have
\begin{equation*}
F^{-1}(1-u)^{\prime }=-1/(u\log u)
\end{equation*}
and let
\begin{equation*}
s(u)=uF^{-1}(1-u)^{\prime }=-1/\log u
\end{equation*}
so that
\begin{equation*}
s^{\prime }(u)=-1/(u(\log u)^{2}).
\end{equation*}
Finally
\begin{equation*}
b(u)=-(\log u)^{-2}\rightarrow 0.
\end{equation*}
With the representation
\begin{equation*}
F^{-1}(1-u)=c+\int_{u}^{1}\frac{s(t)}{t}dt,
\end{equation*}
where $s$ is slowly varying at zero, we get that
\begin{equation*}
\frac{F^{-1}(1-ux)-F^{-1}(1-u)}{s(u)}\rightarrow -\log x.
\end{equation*}
We can use direct methods and get
\begin{equation*}
F^{-1}(1-ux)-F^{-1}(1-u)=\log (\log (ux)/\log u).
\end{equation*}
We remark $v(u,x)=\log (ux)/\log u\rightarrow 1,$ and that $v(u,x)-1=(\log
x)/\log u.$ We may use the expansion of the logarithm function and get
\begin{eqnarray*}
F^{-1}(1-ux)-F^{-1}(1-u) &=&(v-1)-(v-1)^{2}+O\left( (v-1)^{3}\right) \\
&=&(\log x)/\log u+((\log x)/\log u)^{2}/2+O((\log x)/\log u)^{3}).
\end{eqnarray*}
By putting
\begin{equation*}
S(u)=s(u)=1/(\log u),
\end{equation*}
It comes that 
\begin{equation*}
\frac{\dfrac{F^{-1}(1-ux)-F^{-1}(1-u)}{s(u)}+\log x}{S(u)}\rightarrow (\log
x)^{2}/2.
\end{equation*}

\subsubsection{Normal standard}

\noindent Let F be the $d.f.$ \ of a standard normal law. We have the simple
approximation, for $M=\sqrt{2\pi },$ for $x>1,$%
\begin{equation*}
M^{-1}(x^{-1}-x^{-3})e^{-x^{2}/2}\leq 1-F(x)\leq M^{-1}x^{-1}e^{-x^{2}/2}.
\end{equation*}
For $s=1-F(x),$%
\begin{equation*}
-\log M-\log x-\log (1-x^{-2})-x^{2}/2\leq s\leq -\log M-\log x-x^{2}/2
\end{equation*}
\begin{equation*}
-\log M-\log x+\frac{1}{x}+O(x^{-2})-x^{2}/2\leq \log s\leq -\log M-\log
x-x^{2}/2
\end{equation*}
\begin{equation*}
\log M+\log x+x^{2}/2\leq \log (1/s)\leq \log M+\log x-\frac{1}{x}%
+O(x^{-2})+x^{2}/2.
\end{equation*}
And as $x\rightarrow 0\iff s\rightarrow 0,$%
\begin{equation*}
x=F^{-1}(1-s)=(2\log 1/s)^{1/2}(1+o(1)).
\end{equation*}
We easily see that the $o(1)$ term is at least of order $\left( \log
1/s\right) ^{-1}$. This gives
\begin{equation*}
x^{2}/2+\log M+\log x\leq \log (1/s)\leq x^{2}/2+\log M+\log x-\frac{1}{x}%
+O(x^{-2})
\end{equation*}
\begin{equation*}
\log (1/s)+\log (1/s)\times o(1)+\log M+\log ((2\log 1/s)^{1/2}(1+o(1)).
\end{equation*}
The left term is
\begin{equation*}
=\log (1/s)(1+o(1)+\log M(1/2)\log 2+(1/2)\log \log (1/s)+o(1).
\end{equation*}
The right term is
\begin{eqnarray*}
&&\log (1/s)(1+o(1))+\log M+(1/2)\log 2+(1/2)\log \log (1/s) \\
&&\ \ \ \ \ \ \ \ \ \ \ \ \ \ \ \ \ \ \ \ \ \ \ -(2\log
1/s)^{-1/2}(1+o(1))+O((\log 1/s)^{-1})+o(1).
\end{eqnarray*}
The middle term is
\begin{equation*}
(1+o(1))x^{2}/2.
\end{equation*}
By dividing by $\log (1/s)$, we get 
\begin{eqnarray*}
&&1+\frac{(1/2)\log 4\pi +(1/2)\log \log (1/s)+(2\log
1/s)^{-1/2}(1+o(1))+O((\log 1/s)^{-1})}{\log (1/s)} \\
&&\ \ \ \ \ \ \ \ \ \ \ \ \ \ \ \ \ \ \ \ \ \ \ \ \ \ \ 
\begin{tabular}{l}
$\leq $%
\end{tabular}
(x/(2\log (1/s))^{1/2})^{2}+o(1)\times (x/(2\log (1/s))^{1/2})^{2} \\
&&\ \ \ \ \ \ \ \ \ \ \ \ \ \ \ \ \ \ \ \ \ \ \ \ \ \ \ 
\begin{tabular}{l}
$\leq $%
\end{tabular}
1+\frac{(1/2)\log 4\pi +(1/2)\log \log (1/s)}{\log (1/s)}.
\end{eqnarray*}
Then
\begin{eqnarray*}
&&1+\frac{(2\log 1/s)^{-1/2}(1+o(1))+O((\log 1/s)^{-1})+\log (1/s)\times o(1)%
}{\log (1/s)} \\
&&\ \ \ \ \ \ \ \ \ \ \ 
\begin{tabular}{l}
$\leq $%
\end{tabular}
(x/(2\log (1/s))^{1/2})^{2}-\frac{(1/2)\log 4\pi +(1/2)\log \log (1/s)}{\log
(1/s)} \\
&&\ \ \ \ \ \ \ \ \ \ \ 
\begin{tabular}{l}
$=$%
\end{tabular}
(x/(2\log (1/s))^{1/2})^{2}-\frac{(1/2)\log 4\pi +(1/2)\log \log (1/s)}{\log
(1/s)} \\
&&\ \ \ \ \ \ \ \ \ \ \ 
\begin{tabular}{l}
$=$%
\end{tabular}
1+\frac{o(1)}{\log (1/s)}.
\end{eqnarray*}
Then
\begin{equation*}
x^{2}=(2\log (1/s)\left\{ 1-\frac{(1/2)\log 4\pi +(1/2)\log \log (1/s)+o(1)}{%
\log (1/s)}\right\} ,
\end{equation*}
and
\begin{equation*}
x=F^{-1}(1-u)=(2\log (1/s))^{1/2}\left\{ 1-\frac{(1/2)\log 4\pi +(1/2)\log
\log (1/s)+o(1)}{2\log (1/s)}\right\} .
\end{equation*}
We have
\begin{equation*}
F^{-1}(1-s)-F^{-1}(1-xs)=A(x,s)+B(x,s)+C(x,s)
\end{equation*}
with
\begin{equation*}
A(x,s)=(2\log (1/s))^{1/2}-(2\log (x/s))^{1/2}=-(2\log (1/s))^{1/2}((\frac{%
\log (x/s)}{\log (1/s)})^{1/2}-1).
\end{equation*}
But
\begin{eqnarray*}
\left( \frac{\log (x/s)}{\log (1/s)}\right) ^{1/2} &=&\left( 1+\left( \frac{%
\log (x/s)}{\log (1/s)}-1\right) \right) ^{1/2} \\
&=&1+\frac{1}{2}\left( \frac{\log (x/s)}{\log (1/s)}-1\right) -\frac{1}{8}%
\left( \frac{\log (x/s)}{\log (1/s)}-1\right) ^{2}+O\left( \left( \frac{\log
(x/s)}{\log (1/s)}-1\right) ^{3}\right) \\
&=&1+\frac{1}{2}\frac{\log x}{\log 1/s}-\frac{1}{2}\frac{(\log x)^{2}}{%
\left( 2\log 1/s\right) ^{2}}+O\left( \left( \log 1/s\right) ^{-3}\right) .
\end{eqnarray*}
We get
\begin{equation*}
A(x,s)=-\frac{\log x}{\left( 2\log 1/s\right) ^{1/2}}-\frac{1}{2}\frac{(\log
x)^{2}}{\left( 2\log 1/s\right) ^{3/2}}+O\left( (\log 1/s)^{-5/2}\right)
\end{equation*}
and
\begin{eqnarray*}
B(x,s) &=&\frac{(1/2)\log 4\pi +(1/2)\log \log (x/s)+o(\left( \log
1/s\right) ^{-1})}{\left( 2\log (x/s)\right) ^{1/2}} \\
&&-\frac{(1/2)\log 4\pi +(1/2)\log \log (1/s)+o(\left( \log 1/s\right) ^{-1})%
}{\left( 2\log (x/s)\right) ^{1/2}},
\end{eqnarray*}
\begin{eqnarray*}
C(x,s) &=&\frac{(1/2)\log 4\pi +(1/2)\log \log (1/s)+o(\left( \log
1/s\right) ^{-1})}{\left( 2\log (x/s)\right) ^{1/2}} \\
&&-\frac{(1/2)\log 4\pi +(1/2)\log \log (1/s)+o(\left( \log 1/s\right) ^{-1})%
}{\left( 2\log (1/s)\right) ^{1/2}}
\end{eqnarray*}
and
\begin{equation*}
B(x,s)=\frac{1}{\left( 2\log (x/s)\right) ^{1/2}}(\frac{1}{2}\log (\log
(x/s)/\log (1/s))+\left( o(\left( \log 1/s\right) ^{-1})\right) .
\end{equation*}
But
\begin{eqnarray*}
\log (\log (x/s)/\log (1/s)) &=&\log (1+\left( \log (x/s)/\log
(1/s)-1\right) ) \\
&=&\left( \log (x/s)/\log (1/s)-1\right) -\frac{1}{2}\left( \log (x/s)/\log
(1/s)-1\right) ^{2} \\
&&+O\left( \left( \log (x/s)/\log (1/s)-1\right) ^{3}\right) \\
&=&\frac{\log x}{\log 1/s}-\frac{1}{2}\frac{\left( \log x\right) ^{2}}{%
\left( \log 1/s\right) ^{2}}+O((\log 1/s)^{-3}).
\end{eqnarray*}
Then
\begin{eqnarray*}
B(x,s) &=&\frac{1}{(2\log 1/s)^{1/2}}\left\{ \frac{\log x}{\log 1/s}-\frac{1%
}{2}\frac{\left( \log x\right) ^{2}}{\left( \log 1/s\right) ^{2}}+O((\log
1/s)^{-3})\right\} \\
&&+\left( o(\left( \log 1/s\right) ^{-3/2})\right) ,
\end{eqnarray*}
and
\begin{eqnarray*}
C(x,s) &=&\frac{(1/2)\log 4\pi +(1/2)\log \log (1/s)+o(\left( \log
1/s\right) ^{-1})}{\left( 2\log (x/s)\right) ^{1/2}} \\
&&\times \left\{ \left( 2\log (1/s)\right) ^{1/2}-\left( 2\log (x/s)\right)
^{1/2}\right\} \\
&=&\frac{(1/2)\log 4\pi +(1/2)\log \log (1/s)+o(\left( \log 1/s\right) ^{-1})%
}{\left( 2\log (1/s)\right) ^{1/2}} \\
&&\times \left\{ -\frac{\log x}{\left( 2\log 1/s\right) ^{1/2}}-\frac{1}{2}%
\frac{(\log x)^{2}}{\left( 2\log 1/s\right) ^{3/2}}+O((\log
1/s)^{-5/2})\right\} .
\end{eqnarray*}
Recall
\begin{equation*}
A(x,s)=\left\{ -\frac{\log x}{\left( 2\log 1/s\right) ^{1/2}}-\frac{1}{2}%
\frac{(\log x)^{2}}{\left( 2\log 1/s\right) ^{3/2}}+O((\log
1/s)^{-5/2})\right\} .
\end{equation*}
We conclude that 
\begin{equation}
\frac{\left( 2\log (1/s)\right) ^{1/2}}{(1/2)\log 4\pi +(1/2)\log \log (1/s)}%
\left\{ \frac{F^{-1}(1-s)-F^{-1}(1-xs)}{(2\log (1/s))^{-1/2}}+\log x\right\}
\rightarrow -\log x.  \label{norm}
\end{equation}

\subsubsection{Lognormal}

\noindent We have
\begin{equation*}
G^{-1}(1-s)=\exp (F^{-1}(1-s)),\text{ }s\in (0,1),
\end{equation*}
where $F$ is standard normal. This gives
\begin{equation*}
G^{-1}(1-u)-G^{-1}(1-ux)=\exp (F^{-1}(1-s))\times (1-\exp
(F^{-1}(1-xs)-F^{-1}(1-s))).
\end{equation*}
But
\begin{eqnarray*}
\exp (F^{-1}(1-xs)-F^{-1}(1-s)) &=&\exp \left( -\frac{%
F^{-1}(1-s)-F^{-1}(1-xs)}{(2\log (1/s))^{1/2}}(2\log (1/s))^{1/2}\right) \\
&=&\exp \left( -\frac{D(x,s)}{(2\log (1/s))^{1/2}}\right) \\
&=&1-\frac{D(x,s)}{(2\log (1/s))^{1/2}}+\frac{1}{2}\frac{D(x,s)^{2}}{(2\log
(1/s))}+\left( -\frac{1}{(2\log (1/s))^{3/2}}\right) ,
\end{eqnarray*}
where $D(x,s)$ is defined in (\ref{norm}). This yields
\begin{equation*}
\frac{G^{-1}(1-s)-G^{-1}(1-sx)}{\exp (F^{-1}(1-s))(2\log (1/s))^{-1/2}}%
=D(x,s)-\frac{1}{2}\frac{D(x,s)^{2}}{(2\log (1/s))^{1/2}}+O\left( -\frac{1}{%
(2\log (1/s))^{1}}\right) .
\end{equation*}
Thus
\begin{equation*}
\frac{\left( 2\log (1/s)\right) ^{1/2}}{(1/2)\log 4\pi +(1/2)\log \log (1/s)}%
\left\{ \frac{G^{-1}(1-s)-G^{-1}(1-sx)}{\exp (F^{-1}(1-s))(2\log
(1/s))^{-1/2}}+\log x\right\} \rightarrow -\frac{1}{2}\left( \log x\right)
^{2}.
\end{equation*}

\subsubsection{Logistic law}

\begin{equation*}
F(x)=1-2/(1+e^{x}),
\end{equation*}
that is
\begin{equation*}
F^{-1}(1-s)=\log s^{-1}(2-s).
\end{equation*}
Routine computations yield
\begin{eqnarray*}
F^{-1}(1-s)-F^{-1}(1-xs) &=&\log x+\log (2-s)/(2-xs)=\log x+\log
(1+(s(x-1)/(2-xs)) \\
&=&\log x+(s(x-1)/(2-xs))+O(s^{2}).
\end{eqnarray*}
Thus
\begin{equation*}
(s/2)^{-1}\left\{ F^{-1}(1-s)-F^{-1}(1-xs)-\log x\right\} \rightarrow \frac{%
x-1}{2}.
\end{equation*}

\subsubsection{Log-Expo}

\noindent We have 
\begin{equation*}
F^{-1}(1-s)=\log \log (1/s),\text{ }s\in (0,1).
\end{equation*}
\begin{eqnarray*}
F^{-1}(1-s)-F^{-1}(1-xs) &=&\log ((\log (x/s)/(\log (1/s))) \\
&=&\log (1+\left( \log (x/s)/\log (1/s)-1\right) ) \\
&=&\left( \log (x/s)/\log (1/s)-1\right) -\frac{1}{2}\left( \log (x/s)/\log
(1/s)-1\right) ^{2} \\
&&+O(\left( \log (x/s)/\log (1/s)-1\right) ^{3}) \\
&=&\frac{\log x}{\log \left( 1/s\right) }-\frac{1}{2}\frac{\left( \log
x\right) ^{2}}{\left( \log \left( 1/s\right) \right) ^{2}}+O((\log \left(
1/s\right) )^{-3}).
\end{eqnarray*}
This gives
\begin{equation*}
(\log \left( 1/s\right) )^{2}\left\{ \frac{F^{-1}(1-s)-F^{-1}(1-xs)}{(\log
\left( 1/s\right) )^{-1}}-\log x\right\} \rightarrow -\frac{1}{2}(\log x)^{2}
\end{equation*}

\subsubsection{Reversed Burr's df}

\noindent We have 
\begin{equation*}
F(x)=1-((-x)^{-\rho /\gamma }+1)^{1/\rho },\text{ }x\leq 0,\text{ }\rho <0%
\text{ and }\gamma >0.
\end{equation*}
Then
\begin{equation*}
F^{-1}(1-u)=-(u^{\rho }-1)^{-\gamma /\rho }.
\end{equation*}
and
\begin{eqnarray*}
F^{-1}(1-u)-F^{-1}(1-ux) &=&((xu)^{\rho }-1)^{-\gamma /\rho }-(u^{\rho
}-1)^{-\gamma /\rho } \\
&=&(u^{\rho }-1)^{-\gamma /\rho }\left\{ \left( \frac{x^{\rho }u^{\rho }-1}{%
u^{\rho }-1}\right) ^{-\gamma /\rho }-1\right\} \\
&=&(u^{\rho }-1)^{-\gamma /\rho }\left\{ x^{-\gamma }\left( \frac{1-x^{-\rho
}u^{-\rho }}{1-u^{-\rho }}\right) ^{-\gamma /\rho }-1\right\} .
\end{eqnarray*}
But
\begin{eqnarray*}
\left( \frac{1-x^{-\rho }u^{-\rho }}{1-u^{-\rho }}\right) ^{-\gamma /\rho }
&=&\left( 1+\left\{ \frac{1-x^{-\rho }u^{-\rho }}{1-u^{-\rho }}-1\right\}
\right) ^{-\gamma /\rho } \\
&=&\left( 1+\left\{ \frac{(1-x^{-\rho })u^{-\rho }}{1-u^{-\rho }}\right\}
\right) ^{-\gamma /\rho } \\
&=&1-\frac{\gamma }{\rho }\frac{(1-x^{-\rho })u^{-\rho }}{1-u^{-\rho }}%
+O(u^{-2\rho }).
\end{eqnarray*}
Thus
\begin{equation*}
\frac{F^{-1}(1-u)-F^{-1}(1-ux)}{\gamma (u^{\rho }-1)^{-\gamma /\rho }}=\frac{%
x^{-\gamma }-1}{\gamma }-x^{-\gamma }\frac{1}{\rho }\frac{(1-x^{-\rho
})u^{-\rho }}{1-u^{-\rho }}+O(x^{-\gamma }\gamma ^{-1}u^{-2\rho }).
\end{equation*}
So
\begin{equation*}
\left( u^{-\rho }\right) ^{-1}\left\{ \frac{F^{-1}(1-u)-F^{-1}(1-ux)}{\gamma
(u^{\rho }-1)^{-\gamma /\rho }}-\frac{x^{-\gamma }-1}{\gamma }\right\} =-%
\frac{x^{-\gamma }(1-x^{-\rho })}{\rho }
\end{equation*}
We now summarize the results of these computations in the next subsection.

\subsection{Tables of functions $s$ and $b$}

\label{sstable}

\begin{center}
\begin{tabular}{llll}
\textbf{Name\ \ \ \ \ \ \ \ \ \ \ \ \ \ \ \ \ \ \ } & $F$ &  & $b$ \\ 
\hline\hline
Burr & $1-(x^{-\rho /\gamma }+1)^{1/\rho },$ & $x\geq 0,\rho <0,\gamma >0$
& $u^{\rho }$ \\ 
Reversed Burr & $1-((-x)^{-\rho /\gamma }+1)^{1/\rho },$ & $x\leq 0,\rho
<0,\gamma <0$ & $u^{-\rho }$ \\ 
Singh-Maddala & $1-F(x)=(1+ax^{b})^{-c},$ & $x\geq 0$ & $(bc(1-u^{-1/c})^{-1}$ \\ 
Log-Sm & $1-G(x)=(1+ae^{xb})^{-c}$ &  & $\frac{u^{1/c}}{(bc)(1-u^{1/c})}$
\\ 
Exponentiel & $1-e^{-x},$ & $x\geq 0$ & $(\log u)^{-1}$ \\ 
Log-Expo & $1-F(x)=\exp (-e^{-x})$ &  & $1/\log u$ \\ 
Normal & $\phi (x)$ &  & $(\log 1/u)^{-3/2}$ \\ 
Lognormal & $\phi (e^{x})$ &  & $(\log 1/u)^{-1}$ \\ 
Logistic & $1-2/(1+e^{x}),$ & $x\geq 0$ & $(\log 1/u)^{-1}$ \\ \hline\hline
&  &  & 
\end{tabular}
\bigskip\ 
\begin{tabular}{llll}
\textbf{Name\ \ \ \ \ \ \ \ \ \ \ \ \ \ \ \ \ \ \ } & $gamma$ & $h_{\gamma
,\rho }(x)$ & $S(u)$ \\ \hline\hline
Burr & $1/\gamma $ & $-\frac{(x^{\rho }-1)x^{-\gamma -\rho }}{\rho }$ & $%
(u^{\rho }-1)^{-1}$ \\ 
Reversed Burr & $-1/\gamma $ & $-\frac{x^{-\gamma }(1-x^{-\rho })}{\rho }$ & 
$u^{-\rho }$ \\ 
Singh-Maddala & $1/(bc)$ & $\frac{\left( x^{-1/c}-1\right) x^{-1/(bc)+1/c}}{%
1/c}$ & $-(1-u^{-1/c})^{-1}$ \\ 
Log-Sm & $0$ & $-c(x^{-1/c}-1)$ & $cu^{1/c}/(\gamma b)$ \\ 
Exponentiel & $0$ & Not applicable & Not applicable \\ 
Log-Expo & $0$ & $-(\log x)^{2}/2$ & $1/\log u$ \\ 
Normal & $0$ & $-\log x$ & $D(u)=\left\{ \frac{(1/2)\log 4\pi +(1/2)\log
\log (1/u)}{\left( 2\log (1/u)\right) ^{1/2}}\right\} ^{-1}$ \\ 
Lognormal & $0$ & $(\log x)^{2}/2$ & the same $D(u)$ \\ 
Logistic & $0$ & $x-1$ & $u/2$ \\ \hline\hline
&  &  & 
\end{tabular}
\end{center}

Statistical applications

\noindent Pratically, the normality results on statistics based on the
extremes are applied for ultimately differentiable distribution functions
(at $+\infty ).$ They usually depend of the functions $b(.)$\ for in the
representation scheme and, on $S$ in the second order condition one. This
means that we may move from one approach to the other. Let us illustrate
this with two examples.

\subsection{Large quantiles process}

\noindent Let $X_{1},X_{2},$\ ... be a sequence of real and \ independant
random variables indentically distributed and associated to the distribution
function $F(x)=P(X_{i}$\ $\leq x),x\in R.$ We suppose that these random
variables are represented as $X_{i}=F^{-1}(U_{i}),$ $i=1,2,...,$ where $%
U_{1},U_{2},...$ are standard uniform independant random variables. For each 
$n,$ $U_{1,n}<...<U_{n,n}$ denote the order statistics based on $%
U_{1},...,U_{n}$. Finally let $\alpha >0$ and $a>0$ and
\begin{equation*}
k\rightarrow \infty ;\text{ }k/n\rightarrow 0\text{ }and\text{ }\log \log
n/k\rightarrow 0\text{ as }n\rightarrow \infty .
\end{equation*}
Consider this large quantile proccess (see Drees \cite{drees})
\begin{equation*}
A_{n}(s,\alpha )=X_{n-[k/s^{\alpha }]+1,n}-F^{-1}(1-[k/s^{\alpha }]/n).
\end{equation*}
We suppose that $F$ is in the extremal domain. We use first the
representation scheme.

\subsubsection{Representation approach}

\noindent Consider the function $p$ and $b$ defined in Theorem \ref{A}. For
any $\lambda >1,$ put the convention 
\begin{equation*}
d_{n}(h,a,\alpha )=\sup \{\left\vert h(t)\right\vert ,\left\vert
t\right\vert \leq \lbrack a^{-\alpha }]\lambda k/n.
\end{equation*}
We may then define the regularity condition,
\begin{equation}
\sqrt{k}(d_{n}(p,a,\alpha )\vee d_{n}(b,a,\alpha ))\rightarrow 0  \tag{RCREP}
\end{equation}
under which we may find a uniform Gaussian approximation of $A_{n}(s,\alpha
).$ Put for convenience
\begin{equation*}
k(s,\alpha )=[k/s^{\alpha }],\text{ }l(n,\alpha )=k(s,\alpha )/n,\ \
k^{\prime }(s,\alpha )=k(s,\alpha )/k,\text{ and }U_{[k/s^{\alpha
}],n}=U_{k(s,\alpha ),n}.
\end{equation*}
For 
\begin{equation*}
a(k/n)=c(1+p(k/n))(k/n)^{-\gamma }\exp (\int_{k/n}^{1}b(t)dt),
\end{equation*}
we have
\begin{eqnarray*}
A_{n}(1) &=&X_{n-[k/s^{\alpha }]+1,n}/a(k/n) \\
&=&(1+O(d_{n}(p,\alpha ,\lambda ))(1+O(d_{n}(p,\alpha ,\lambda ))\times
(U_{k(s,\alpha ),n}/l(n,s,\alpha )^{-\gamma }k^{\prime }(s,\alpha )^{-\gamma
}.
\end{eqnarray*}
We also have

\begin{eqnarray*}
A_{n}(2) &=&F^{-1}(1-l(n,s,\alpha ))/a(k/n) \\
&=&(1+O(d_{n}(p,\alpha ,\lambda ))(1+O(d_{n}(b,\alpha ,\lambda ))k^{\prime
}(s,\alpha )^{-\gamma }.
\end{eqnarray*}
It follows, since $k^{\prime }(s,\alpha )^{-\gamma }=s^{\alpha \gamma
}(1+O(k^{-1}))$ uniformly in $s\in (a,1)$, that 
\begin{eqnarray*}
A_{n}(1)-A_{n}(2) &=&s^{\alpha \gamma }\left\{ (l(n,s,\alpha
)^{-1}U_{k(s,\alpha ),n})^{-\gamma }-1\right\} \\
&&+O(d_{n}(p,\alpha ,\lambda )\vee d_{n}(b,\alpha ,\lambda ))+O(k^{-1});
\end{eqnarray*}
and by (\cite{adja})
\begin{eqnarray*}
\sqrt{k}\left\{ A_{n}(1)-A_{n}(2)\right\} &=&-\gamma s^{\alpha \gamma }\sqrt{%
k}\left\{ (l(n,s,\alpha )^{-1}U_{k(s,\alpha ),n})-1\right\} \\
&&+O(\sqrt{k}d_{n}(p,\alpha ,\lambda ))+O(\sqrt{k}d_{n}(b,\alpha ,\lambda )))
\\
&=&-\gamma s^{\alpha \gamma }W_{n}(1,s^{\alpha })+o_{p}(a,\alpha ,s)
\end{eqnarray*}
whenever (RCREP) is valid. We then obtain the limiting law of the process of
large quantiles under this condition. When $F$ is differentiable in the
neighborhood of $+\infty ,$ we may take $p=0$ and (RCREP) becomes
\begin{equation*}
\sqrt{k}d_{n}(b,a,\alpha )\rightarrow 0.
\end{equation*}

\noindent Under this (RCREP), the large quantile process behaves as the
Gaussian stochastic process $-\gamma s^{\alpha \gamma }(s,\alpha
)W_{n}(1,s^{\alpha })$.

\subsubsection{Second order condition approach}

\noindent There exist functions $a(\cdot )$ and $S(\cdot )$ ($a(\cdots )$ is
not necessarily the same as the previous function $s(\cdots )$), such that
the SOC holds. But for statistical purposes, it is more convenient to use
the continuous second order condition, that is for $u_{n}\rightarrow 0,$ for 
$x_{n}\rightarrow x>0,$%
\begin{equation*}
\lim_{n\rightarrow \infty }S(u_{n})\left\{ \frac{%
F^{-1}(1-u_{n})-F^{-1}(1-u_{n}x_{n})}{a(u_{n})}-\frac{x_{n}^{-\gamma }-1}{%
\gamma }\right\} =h_{\gamma }(x).
\end{equation*}
A simple argument based on compactness yields for $u_{n}\rightarrow 0$ and
for $0<a<b$,
\begin{equation*}
\lim_{n\rightarrow \infty }\sup_{a\leq x\leq b}\left\vert S(u_{n})\left\{ 
\frac{F^{-1}(1-u_{n})-F^{-1}(1-u_{n}x)}{a(u_{n})}-\frac{x^{-\gamma }-1}{%
\gamma }\right\} -h_{\gamma }(x)\right\vert =0.
\end{equation*}
Put $x(n,s,\alpha )=l(n,s,\alpha )^{-1}U_{k(s,\alpha ),n}\rightarrow 1.$ For 
$s\geq a,$ we may see that $l(n,s,\alpha )^{-1}U_{k(s,\alpha
),n}=1+k^{-1/2}(W_{n}(s^{\alpha })+o_{P}(1))$, uniformly in $s\in (a,1)$,
where $W_{n}$ is a standard Wiener process (see Lemma 1 in \cite{adja}).
Then we may apply the CSOC as follows : 
\begin{eqnarray*}
&&\sup_{a\leq s\leq b}S(l(n,s,\alpha ))^{-1}\left\vert \left\{ \frac{%
F^{-1}(1-l(n,s,\alpha ))-F^{-1}(1-U_{k(s,\alpha ),n})}{a(l(n,s,\alpha ))}%
\right. \right. \\
&&\ \ \ \ \ \ \ \ \ \ \ \ \ \ \ \ \ \ \ \ \ \ \ \ \ \ \ -\left. \left. \frac{%
(l(n,s,\alpha )^{-1}U_{k(s,\alpha ),n})^{-\gamma }-1}{\gamma }\right\}
-h_{\gamma }(x(n,s,\alpha ))\right\vert \\
&&\ \ \ \ \ \ \ \ \ \ \ \ \ \ \ 
\begin{tabular}{l}
$=$%
\end{tabular}
0
\end{eqnarray*}

\noindent This gives, uniformly in $s\in (a,1),$%
\begin{eqnarray*}
\frac{F^{-1}(1-l(n,s,\alpha ))-F^{-1}(1-U_{k(s,\alpha ),n})}{a(l(n,s,\alpha
))} &=&\frac{(l(n,s,\alpha )^{-1}U_{k(s,\alpha ),n})^{-\gamma }-1}{\gamma }
\\
&&-(h_{\gamma }(x(n,s))+o_{P}(1))S(l(n,s,\alpha )).
\end{eqnarray*}
Then
\begin{eqnarray*}
\frac{\sqrt{k}\left\{ F^{-1}(1-l(n,s))-F^{-1}(1-U_{k(s,\alpha ),n})\right\} 
}{a(l(n,s,\alpha ))} &=&\frac{\sqrt{k}(l(n,s)^{-1}U_{k(s,\alpha
),n})^{-\gamma }-1}{\gamma } \\
&&-(h_{\gamma }(1)+o_{P}(1))S(l(n,s))\sqrt{k}.
\end{eqnarray*}
We will apply Lemma 1\ in (\cite{adja}). Since $S(l(n,s,\alpha
))=O(S(l(n,s)) $ and $a(l(n,s,\alpha ))\sim s^{\alpha \gamma }a(l(n,1,1)),$
we also get
\begin{equation*}
\frac{\sqrt{k}\left\{ F^{-1}(1-l(n,s))-F^{-1}(1-U_{k(s,\alpha ),n})\right\} 
}{a(l(n,s,\alpha ))}
\end{equation*}
\begin{equation*}
=-W_{n}(1,s^{\alpha })-(h_{\gamma }(1)+o_{P}(1))S(l(n,s,\alpha ))\sqrt{k}.
\end{equation*}
\begin{equation*}
\frac{\sqrt{k}\left\{ F^{-1}(1-l(n,s))-F^{-1}(1-U_{k(s,\alpha ),n})\right\} 
}{a(l(n,1,1))}=-s^{\alpha \gamma }W_{n}(1,s^{\alpha })+(h_{\gamma
}(1)+o_{P}(1))\times S(k/n))\sqrt{k}.
\end{equation*}
We get the regularity condition 
\begin{equation}
\sqrt{k}S(k/n)\rightarrow 0.  \tag{RCSOC}
\end{equation}

\begin{conclusion}
In both cases, we conclude that the large quantile process behaves as the
Gaussian process $-s^{\alpha \gamma }W_{n}(1,s^{\alpha })$ when
appropriately normalized under conditions based on $b$ or on $S$.

\noindent By comparing (RCREP) and (RCSOC), we see that the present
normality result in the representation scheme uses the function $b$ while
the Second order one relies on $S$. In fact, almost all the normality
results in both cases rely either on $b$ in the Representation scheme or on $%
S$ in the Second order model. We also see that the second order scheme seems
to use a shorter way. But, as a compensation, the function $S$, as we may
see it here, is more complicated to get. Indeed for differentiable
distribution functions, the function $b$, is easiliy obtained.
\end{conclusion}

\subsection{Functional Hill process}

\subsubsection{Representation approach}

\noindent Now consider the functional Hill process 
\begin{equation*}
T_{n}(f)=\sum_{j=1}^{j=k}f(j)\left( \log X_{n-j+1,n}-\log X_{n-j,n}\right) ,
\end{equation*}
where f is some positive and bounded function and $k=k(n)$ is a sequence if
positive integer such that $1\leq 1\leq k\leq n$ and $k/n\rightarrow 0$ as $%
n\rightarrow +\infty $. We are going to study the process under the
hypothesis $F\in D(\psi _{1/\gamma })=D(G_{-\gamma }),$ $\gamma >0$. Now
using the same representation $X_{i}=F^{-1}(U_{i}),$ $i=1,2,...$ We get 
\begin{eqnarray*}
T_{n}(f) &=&\sum_{j=1}^{k}f(j)(\log X_{n-j+1,n}-\log
X_{n-j,n})=\sum_{j=1}^{k}f(j)(-(y_{0}-\log X_{n-j+1,n})+(y_{0}-\log
X_{n-j,n})) \\
&=&\sum_{j=1}^{k}f(j)\left\{ c(1+p(U_{j+1,n}))(U_{j+1,n})^{\gamma }\exp
\left( \int_{U_{j+1,n}}^{1}\frac{b(t)}{t}dt\right) \right. \\
&&\ \ \ \ \ \ \ \ \ \ \ \ \ \ \ \ \ \ \ \ \ \ \ \ \ \ \ \ \ \ \ \ \ \left.
-c(1+p(U_{j,n}))(U_{j,n})^{\gamma }\exp \left( \int_{U_{j,n}}^{1}\frac{b(t)}{%
t}dt\right) \right\} \\
&=&\sum_{j=1}^{k}f(j)\left\{ c(1+p(U_{j,n}))(U_{j,n})^{1/\gamma }\exp \left(
\int_{U_{j,n}}^{1}\frac{b(t)}{t}dt\right) \right\} \\
&&\ \ \ \ \ \ \ \ \ \ \ \ \ \ \times \left\{ \frac{1+p(U_{j+1,n})}{%
1+p(U_{j,n})}\left( \frac{U_{j+1,n}}{U_{j,n}}\right) ^{\gamma }\exp \left(
\int_{U_{j+1,n}}^{U_{j,n}}\frac{b(t)}{t}dt\right) -1\right\} .
\end{eqnarray*}

\noindent But 
\begin{equation*}
\left( \frac{U_{j+1,n}}{U_{j,n}}\right) ^{\gamma }=\exp \left( \frac{\gamma 
}{j}\log \left( \frac{U_{j+1,n}}{U_{j,n}}\right) ^{j}\right) =\exp \left( 
\frac{\gamma }{j}E_{j}\right) =\exp (F_{j})
\end{equation*}
where, by the Malmquist representation (see \cite{shwell}, p. 336), the $%
E_{j}^{\prime }s$ are independent standard exponential random variables. Let
also 
\begin{equation*}
p_{n}=\sup \{\left\vert p(u)\right\vert ,0\leq u\leq U_{k+1,n}\}\rightarrow
_{P}0,
\end{equation*}
\begin{equation*}
b_{n}=\sup \{\left\vert b(u)\right\vert ,0\leq u\leq U_{k+1,n}\}\rightarrow
_{P}0,
\end{equation*}
as $n\rightarrow +\infty $, and $c_{n}=a_{n}\vee (b_{n}\log k)$. Then
\begin{eqnarray*}
\left\{ \frac{1+p(U_{j+1,n})}{1+p(U_{j,n})}\left( \frac{U_{j+1,n}}{U_{j,n}}%
\right) ^{\gamma }\exp \left( \int_{U_{j+1,n}}^{U_{j,n}}\frac{b(t)}{t}%
dt\right) -1\right\} &=&\exp (F_{j})(1+O(p_{n}))\exp (O(b_{n})E_{j}/j)-1 \\
&=&F_{j}(1+O(p_{n}))(1+O(b_{n}\log k)) \\
&=&F_{j}(1+O(c_{n}))-1
\end{eqnarray*}
Let also 
\begin{equation*}
s_{n}=y_{0}-G^{-1}(1-U_{k,n})=y_{0}-\log X_{n-k+1,n}=c(U_{k+1,n})^{1/\gamma
}\left( 1+\exp \left( \int_{U_{j,n}}^{1}\frac{b(t)}{t}dt\right) \right) .
\end{equation*}
This gives 
\begin{equation*}
=T_{n}(f)/s_{n}
\end{equation*}
\begin{equation*}
=\sum_{j=1}^{k}f(j)\left\{ \frac{1+p(U_{j,n})}{1+p(U_{k,n})}\left( \frac{%
U_{j,n}}{U_{k,n}}\right) ^{\gamma }\exp \left( \int_{U_{j,n}}^{U_{k,n}}\frac{%
b(t)}{t}dt\right) \right\} \times \left\{ \exp (F_{j})(1+O(c_{n}))-1\right\}
.
\end{equation*}
Let us remark that 
\begin{equation*}
\log \left( \frac{U_{k,n}}{U_{j,n}}\right) =\log \left(
\prod\limits_{h=j}^{k-1}\frac{U_{h+1,n}}{U_{h,n}}\right) =\sum_{h=j}^{k-1}%
\frac{1}{h}\log \left( \frac{U_{h+1,n}}{U_{h,n}}\right) ^{h}=\sum_{h=j}^{k-1}%
\frac{1}{h}E_{h}=O(\log k)
\end{equation*}
and 
\begin{equation*}
\left( \frac{U_{j,n}}{U_{k,n}}\right) ^{\gamma }=\exp \left( -\gamma
\sum_{h=j}^{k-1}\frac{1}{h}E_{h}\right) =\exp \left(
-\sum_{h=j}^{k-1}F_{h}\right) =F_{j}^{\ast },
\end{equation*}
where $F_{k}^{\ast }=0$. Then 
\begin{eqnarray*}
T_{n}(f)/s_{n} &=&\sum_{j=1}^{k}f(j)\left\{ (1+O(p_{n}))(1+O(b_{n}\log
k))\times \exp \left( -\sum_{h=j}^{k-1}F_{h}\right) \right\} \\
&&\ \ \ \ \ \ \ \ \ \ \ \ \ \ \ \ \times \left\{ \exp (F_{j})-1+O(c_{n})\exp
(F_{j})\right\} \\
&=&\sum_{j=1}^{k}f(j)F_{j}^{\ast }\left( \exp (F_{j})-1\right)
+O(c_{n}F_{j}^{\ast }(\exp (F_{j})-1)) \\
&&\ \ \ \ \ \ \ \ \ \ \ \ \ \ +O(c_{n}^{2}F_{j}^{\ast }(\exp
(F_{j})-1))+O(c_{n}\exp (F_{j}))
\end{eqnarray*}
We conclude that $T_{n}(f)/s_{n}$ behaves as that of $%
\sum_{j=1}^{k}f(j)F_{j}^{\ast }\left( \exp (F_{j})-1\right) $ under
regularity conditions based on the functions $p$ and $b$.

\subsubsection{Second order condition approach}

\noindent Let use the continuous second order condition:
\begin{equation*}
S^{-1}(u_{n})\left\{ \frac{G^{-1}(1-x_{n}u_{n})-G^{-1}(1-u_{n})}{s(u_{n})}-%
\frac{x_{n}^{-\gamma }-1}{\gamma }\right\} =h_{\gamma }(x)+o(1),
\end{equation*}
where x$_{n}\rightarrow 0$ et $u_{n}\rightarrow 0$ as $n\rightarrow \infty $
and $u_{n}=\gamma \left\{ y_{0}(G)-G^{-1}(1-u_{n})\right\} $. We get, for $%
G(x)=F(e^{x})$, $x\in R$,

\begin{equation*}
T_{n}(f)=\sum_{j=1}^{k}f(j)(\log X_{n-j+1,n}-\log X_{n-j,n})
\end{equation*}
\begin{equation*}
=\sum_{j=1}^{k}f(j)\left\{ G^{-1}(1-U_{j,n})-G^{-1}(1-U_{j+1,n})\right\} .
\end{equation*}
Let $u_{n}(j)=U_{j,n}$ et $x_{n}(j)=U_{j+1,n}/U_{j,n}=\exp (F_{j})$. Then 
\begin{eqnarray*}
&&\frac{G^{-1}(1-U_{j,n})-G^{-1}(1-U_{j+1,n})}{s(u_{n}(j))}
\begin{tabular}{l}
$=$%
\end{tabular}
-\frac{G^{-1}(1-x_{n}(j)U_{j,n})-G^{-1}(1-U_{j,n})}{s(u_{n}(j))} \\
&&\ \ \ \ \ \ \ \ \ -\left\{ \frac{%
G^{-1}(1-x_{n}(j)U_{j,n})-G^{-1}(1-U_{j,n})}{s(u_{n}(j))}-\frac{\exp \left( -%
\frac{\gamma }{j}E_{j}\right) -1}{\gamma }\right\} -\frac{\exp \left( -\frac{%
\gamma }{j}E_{j}\right) -1}{\gamma } \\
&&\ \ \ \ \ \ \ \ \ -S(u_{n}(j))h_{\gamma }(\exp (F_{j}))-\frac{%
E_{j}^{-\gamma /j}-1}{\gamma }+o_{p}(A(u_{n}(j)))
\end{eqnarray*}
Let us use 
\begin{eqnarray*}
\frac{T_{n}(f)}{s(u_{n}(k))} &=&\sum_{j=1}^{k}f(j)\frac{s(u_{n}(j))}{%
s(u_{n}(k))}\frac{G^{-1}(1-U_{j,n})-G^{-1}(1-U_{j+1,n})}{s(u_{n}(j))} \\
&=&\sum_{j=1}^{k}f(j)\frac{s(u_{n}(j))}{s(u_{n}(k))}\left\{
-S(u_{n}(j))h_{\gamma }(\exp (F_{j}))-\frac{\exp \left( -\frac{\gamma }{j}%
E_{j}\right) -1}{\gamma }+o_{p}(S(u_{n}(j)))\right\} .
\end{eqnarray*}
Let us apply 
\begin{equation*}
\frac{s(u_{n}(j))}{s(u_{n}(k))}=\left\{ (1+O(p_{n}))(1+O(b_{n}\log k))\times
\exp \left( -\sum_{h=j}^{k-1}F_{h}\right) \right\} .
\end{equation*}
We arrive at 
\begin{eqnarray*}
\frac{\gamma T_{n}(f)}{s(u_{n}(k))} &=&\sum_{j=1}^{k}f(j)(1+O(c_{n}))F_{j}^{%
\ast }\left\{ -S(u_{n}(j))h_{\gamma }(\exp (F_{j}))-(\exp
(F_{j})-1)+o_{p}(A(u_{n}(j)))\right\} \\
&=&\sum_{j=1}^{k}f(j)F_{j}^{\ast }(\exp
(F_{j})-1)+\sum_{j=1}^{k}f(j)O(c_{n}F_{j}^{\ast }(\exp (F_{j})-1)) \\
&&+\sum_{j=1}^{k}f(j)(1+O(c_{n}))F_{j}^{\ast }\left\{ -S(u_{n}(j))h_{\gamma
}(\exp (F_{j}))+o_{p}(S(u_{n}(j)))\right\} .
\end{eqnarray*}

\begin{conclusion}
In both cases, we see that when properly normalized, $T_{n}(f)$ behaves as $%
\sum_{j=1}^{k}f(j)F_{j}^{\ast }\left( \exp (F_{j})-1\right) $ under
regularity conditions based on $p$, $b$ or $S$.

\noindent As for the first example, the SOC approach seems shorter. But here
this latter approach still needs the first one.
\end{conclusion}

\section{Conclusion}

\noindent As a general conclusion, we say : \medskip

\begin{enumerate}
\item  The representation approach is more general.

\item  The second order condition seems to be shorter and more unified.

\item  The computation of $b$ is less complicated than that of $S$.

\item  The representation approach is still used within the second order
approach.

\item  The two approaches may be simultanuously used. \medskip
\end{enumerate}

\noindent We conclude that the two approaches are equivalent and we have
proposed for both cases the computation of $b$ and $S$ for usual
distribution functions. \medskip

\noindent \textbf{Acknowledgement.} The paper was finalized while the first
author was visiting MAPMO, University of Orl\'{e}ans, France, in 2009 and
2011. He expresses his warm thanks to responsibles of MAPMO for kind
hospitality. The second author was partially granted by the project AIRE-SUD
of the Institut Regional de Development (IRD) at the Universit\'{e} Gaston
Berger de Saint-Louis.

\end{document}